\documentclass[11pt]{article}
\usepackage[utf8]{inputenc}
\usepackage{amsfonts}
\usepackage{amsmath}
\usepackage{xr-hyper}
\usepackage{soul}
\usepackage{hyperref}
\hypersetup{colorlinks,%
citecolor= blue,%
filecolor= blue,%
linkcolor=blue,%
urlcolor= blue%
}
\usepackage[normalem]{ulem}
\usepackage{graphicx}
\usepackage[dvipsnames]{xcolor}
\usepackage{subcaption}
\usepackage{fullpage,setspace,authblk,lineno}

 \definecolor{burgundy}{rgb}{0.59, 0.0, 0.09}

\makeatletter
\newcommand*{\addFileDependency}[1]{
  \typeout{(#1)}
  \@addtofilelist{#1}
  \IfFileExists{#1}{}{\typeout{No file #1.}}
}
\makeatother

\graphicspath{{./plots/}}

\title{A fourth-order exponential time differencing scheme with real and distinct poles rational approximation for solving  non-linear reaction-diffusion systems} 

\author[1]{W. K. Attipoe}
\author[2]{A. Kleefeld}
\author[1]{E. O. Asante-Asamani \thanks{easantea@clarkson.edu }}

\affil[1]{Department of Mathematics, Clarkson University, Potsdam NY  13699}
\affil[2]{J\"ulich Supercomputing Centre, Forschungszentrum J\"ulich GmbH\\ 52425 J\"ulich, Germany \& Faculty of Medical Engineering and Technomathematics, University of Applied Sciences Aachen, 52428 J\"ulich, Germany.}

\date{}
\begin{document}
\maketitle
\doublespacing
\begin{abstract}
\noindent A fourth-order, L-stable, exponential time differencing Runge–Kutta type scheme is developed to solve non-linear systems of reaction–diffusion equations with non-smooth data. The new scheme, ETDRK4RDP, is constructed by approximating the matrix exponentials in the ETDRK4 scheme with a fourth-order,  L-acceptable, non-Pad\'e rational function having real and distinct poles. Using real and distinct poles rational functions to construct the scheme ensures efficient damping of spurious oscillations arising from non-smooth initial and/or boundary conditions and a straightforward parallelization. We verify empirically that the new ETDRK4RDP scheme is fourth-order accurate for several reaction diffusion systems with Dirichlet and Neumann boundary conditions and show it to be more efficient than competing exponential time differencing schemes, especially when implemented in parallel, with up to six times speed-up in CPU time.
\end{abstract}

\section{Introduction}
Time-dependent reaction-diffusion equations (RDEs) are mathematical models that describe the spatio-temporal dynamics of substances as they react and diffuse in a fluid medium. They are used as classical models of biological pattern formation \cite{leveque2007finite, tyson1999model}, following pioneering work by Alan Turing.  With the inclusion of an advection term, RDEs have been used to model the movement of pollutants through surface and ground water \cite{kamel2007chemical,zhang2022fourier} as well as monitoring the onset and progression of tumors \cite{anderson2000mathematical, chaplain1993model, gatenby1996reaction}. 

Mathematically, RDEs are described by the system of partial differential equations of the form  \begin{eqnarray}\label{eqn1.1}
    \frac{\partial u(x, t)}{\partial t} &=& D \Delta u(x,t) + f(u(x, t), t), \quad \text{in} \quad \Omega_T \\ 
     u(x, t) &=& g(x, t), \quad \text{on} \quad \Gamma_T  \nonumber   
\end{eqnarray} where $\Omega$ describes a bounded open subset of $\mathbb{R}^d$, $d \in \{1,2,3\}$ and $ \Omega_T := \Omega \times (0, T]$ is a cylindrical domain with parabolic boundary $\Gamma_T :=  \overline{\Omega_T} \backslash\Omega_T $. The solution to the system $u: \overline{\Omega_T} \longrightarrow \mathbb{R}^s$, describes the spatio-temporal change in the concentration of $s-$species, defined to be $g(x,t)$ on the parabolic boundary and interacting through the nonlinear reaction function $f: \Omega_T \longrightarrow \mathbb{R}^s $, assumed to be sufficiently smooth. The species also diffuse at rates determined by the coefficient matrix $D = \mathrm{diag}(d_1, d_2, ..., d_s), \ d_i > 0 \quad \forall i $.

A classical approach for solving Eq.~\eqref{eqn1.1} uses the method of lines, where the Laplacian is initially discretized with  $m^d$ points to obtain the system of ordinary differential equations
 \begin{eqnarray}\label{eqn1.2}
    \frac{dU}{dt} + AU &=& F(U(t), t), \quad U \in \mathbb{R}^{s\times m^d} \\ 
     U(x, 0) &=& U_0(x) \nonumber,      
\end{eqnarray}
which can then be solved with an appropriate time-discretization scheme. Here, the matrix $A \in \mathbb{R}^{(s \cdot m^d)\times(s \cdot m^d)}$ represents the spatial discretization of the diffusion operator $-D \Delta$,  $U(t) = (u_1(t), u_2(t), ..., u_s(t), ..., u_{s \times m^d}(t))^\top$ is the time-dependent solution vector and $F \in \mathbb{R}^{s \cdot m^d}$ a nonlinear map approximating $f$ on the spatial grid.


There are a host of time-stepping methods that solve the resulting stiff ODE systems in Eq.~\eqref{eqn1.2}. Popular among them are the Backward Euler, IMEX, and Explicit Runge-Kutta methods. These methods have various challenges that can make their computation expensive. They may be fully implicit, which requires an iterative solver at each time step; multistep, which requires a history of solutions which must be computed accurately or are explicit and require relatively small time steps to improve stability. The family of exponential time differencing (ETD) schemes seeks to address some of these challenges. ETD schemes are developed from the exact solution to Eq.~\eqref{eqn1.2}, given by the integral 
 \begin{equation}\label{eqn1.3}
    U(t) = \mathrm{e}^{A(t - t_0)}U_0 + \int_{t_0}^t \mathrm{e}^{A(t - \tau)} F(U(\tau), \tau) \ \mathrm{d} \tau, \quad t \in [t_0, T],
\end{equation}
obtained through Duhamel's principle, by first approximating the integral. The significance of these schemes is in the accurate treatment of the linear (diffusion) term via matrix exponential leading to stable evolution of stiff problems \cite{asante2025fourth}. The challenge is in integrating the nonlinear function $F$ which depends on the unknown solution $U$. In 2004, Cox and Matthews introduced the class of ETD Runge-Kutta (ETDRK) schemes \cite{mathews2004numerical} which approximate the integral in Eq.~\eqref{eqn1.3} by first interpolating the non-linear function $F$ with a high-order polynomials and then integrating the resulting expression. The resulting scheme, first developed for diagonalizable matrices $A$ and later extended to general invertible matrices by Kassam and Trefethen \cite{kassam2005fourth}, still suffered from numerical instability and inefficiencies resulting from the computation of the matrix exponential and high powers of matrix inverses.  

A subclass of ETDRK schemes, developed to improve on the computational efficiency of ETDRK schemes use rational functions to approximate the matrix exponentials and eliminate terms involving high powers of matrices through algebraic simplification. A number of these schemes use the standard Pad\'e rational approximations (\cite{kleefeld2012etd} \cite{yousuf2012numerical} \cite{yousuf2009efficient} \cite{khaliq2009smoothing}) while others use non-Pad\'e rational functions having real and distinct poles (RDP) \cite{asante2016real}. Desirable rational approximations to the exponential $\mathrm{e}^{-z}$ should have modulus bounded by 1 (A-acceptable) and approach zero asymptotically as $|z|\rightarrow \infty$ (L-acceptable). The class of Pad\'e rational functions that satisfy these two properties tends to have complex poles which makes their implementation on lower-level programming systems such as Fortran or C++ challenging. The attraction of the class of RDP rational functions is that they are L-acceptable, have real poles and lead to an easily parallelizable algorithm. 

In 2016, Asante-Asamani, Khaliq and Wade \cite{asante2016real} developed the first class of ETD Runge-Kutta schemes using an RDP rational function to approximate the matrix exponential. The scheme is however of second order, limiting its application to RDEs discretized with higher order spatial discretization schemes such as spectral methods. Asante-Asamani, Kleefeld and Wade in 2025 \cite{asante2025fourth} also developed an exponential time differencing scheme with dimensional splitting that uses a Pad\'e(2,2) type rational function to approximate the matrix exponentials. In this work, we develop a fourth-order ETDRK4 scheme referred to as ETDRK4RDP, which uses a fourth-order RDP rational function to approximate the matrix exponentials. Empirical convergence analysis supports the fourth-order accuracy of the scheme across various scalar and multidimensional RDEs with Dirichlet and Neumann boundary conditions. Comparison with the competing L-stable, fourth-order ETDRK4P04 scheme, which employs the Pad\'e (0,4) rational function, shows improved efficiency for multidimensional systems of RDEs and a significant speed up in CPU time when implemented in parallel. 

In Section \ref{section2} we present a fourth-order finite difference scheme used to discretize the Laplacian, in order to obtain the diffusion matrix $A$. In Section \ref{section3} we introduce the RDP rational function and its properties, and present the derivation of the ETDRK4RDP scheme. In Section \ref{section4}, we show and discuss the numerical experiments to validate the convergence and efficiency of serial and parallel versions of the scheme on various test problems. We share some concluding remarks and future work in Section \ref{section5}.  

\section{Spatial Discretization} \label{section2}
We begin by discretizing the second-order spatial derivatives in 1D and extend the results using the Kronecker product of matrices to problems in 2D. To do this, consider Eq.~\eqref{eqn1.1} and discretize the Laplacian $\Delta$ on the domain $[a, b]$. Define a uniform mesh of size $\displaystyle h = \frac{b-a}{m+1}$ and divide the domain into $m + 2$ points with $m\geq 3$, with each point on the grid set to $x_j = a + jh$ with $j = 0, 1, 2, ..., m+1$. To ensure the overall scheme is fourth-order, we approximate the second-order partial derivatives using a standard fourth-order central difference scheme (\cite{mathews2004numerical}, \ p. 339). By writing $w_j = w(x_j)$, for a given function $w(x), w: [a, b] \longrightarrow \mathbb{R}$, we can approximate its second-order derivatives with respect to $x$ at $x_j$ as 
\begin{equation}\label{eqn1.411}
       w^{''}\vert_{x_j} = \frac{1}{12h^2} (-w_{j-2} + 16w_{j-1} - 30w_{j} + 16w_{j+1} - w_{j+2}) + \mathcal{O}(h^4)
\end{equation}
\noindent
$j = 2,3,\ldots,m-1$.\\
The approximation used at the nodes $j = 0, 1, m$ and $m+1$ depends on the boundary conditions imposed. For homogeneous Dirichlet boundary conditions, since the values at the endpoints $x_0$ and $x_{m+1}$ are known, we extrapolate the values $w_{-1}$ and $w_{m+2}$ using a fourth-degree Lagrange interpolating polynomial centered at $x_0$ and $x_{m+1}$. The formulas for the derivatives at $x_1$ and $x_m$ can be summarized as

\begin{equation}\label{eqn1.4}
\begin{split}
    w^{''}|_{x_1} &= \frac{1}{12h^2} (11w_{0} - 20w_{1} + 6w_{2} + 4w_{3} - w_{4})\\
    w^{''}|_{x_m} &= \frac{1}{12h^2} (-w_{m-3} + 4w_{m-2} + 6w_{m-1} - 20w_{m} + 11w_{m+1}).
\end{split}
\end{equation}

\noindent
For homogeneous Neumann boundary conditions, the value of the function at $w_{-1}$ and $w_{m+1}$ is unknown hence we consider Eq.~\eqref{eqn1.411} and set $w_{-1} = w_{1}$ and $w_m = w_{m+2}$ from the boundary condition

$$ \frac{w_1 - w_{-1}}{2h} = 0  \implies w_1 = w_{-1} \quad \text{and} \quad \frac{w_m - w_{m-2}}{2h} + 0  \implies w_m = w_{m+2}. $$
Thus, the approximation at the boundaries summarizes to 
\begin{equation}\label{eqn1.5}
\begin{split}
     w^{''}|_{x_1} &\approx \frac{1}{12h^2} (16w_{0} - 31w_{1} + 16w_{2} - w_{3}),\\
    w^{''}|_{x_m} &\approx \frac{1}{12h^2} (-w_{m-2} + 16w_{m-1} - 30w_{m} = 16w_{m=1} - w_{m+2}).
\end{split}
\end{equation}

\noindent The derivatives at $x_0$ and $x_{m+1}$ are set to 
\begin{equation}\label{eqn1.6}
\begin{split}
    w^{''}|_{x_0} &\approx \frac{1}{12h^2} (-30w_{0} +32w_{1} - 2w_{2}),\\
    w^{''}|_{x_{m+1}} &\approx \frac{1}{12h^2} (-2w_{m-1} + 32w_{m} - 30w_{m+1}).
\end{split}
\end{equation}
\noindent
A similar discretization has been used in \cite{gibou2005fourth} to achieve fourth-order accuracy on irregular domains. Let $I_p$ be a $ p$-dimensional identity matrix and $B_p$ the matrix obtained from the 1D discretization of the Laplacian using Eq.~\eqref{eqn1.411}. For systems in 2D, we make use of the Kronecker product for the discrete Laplacian presented in \cite{hundsdorfer2013numerical} to formulate the system matrix. In particular, for discrete Laplacian $A$ in 2D, we have that $A=A_1 + A_2$, where  
\begin{equation}\label{eqn2.2221}
 A_1 = B_p \otimes I_p  \quad \text{and} \quad A_2 = I_p \otimes B_p  
\end{equation}
with $B_p$ and $I_p$ as the one-dimensional finite difference discretization of the second derivative and identity matrix of equivalent dimension $p$, respectively. For Dirichlet boundary conditions and $s-$species, $p = s \cdot m$, whereas $p = s \cdot (m + 2)$ for Neumann boundary conditions. 

\section{Development of ETDRK4RDP scheme}
\label{section3}
\subsection{Background}
Cox–Matthews \cite{mathews2004numerical} and Kassam–Trefethen \cite{kassam2005fourth} in 2002 developed a class of fourth-order ETD Runge-Kutta (ETDRK4) schemes which approximate the integral in the exact solution Eq.~\eqref{eqn1.3} by interpolating the non-linear function $F$ with higher-order polynomials. The scheme is most preferred due to the accurate treatment of the linear (diffusion) term via matrix exponential leading to stable evolution of stiff problems. The approximation is given by the recurrence relation,

\begin{eqnarray}\label{eqn2.21}
U_{n+1} &=& \mathrm{e}^{-kA}U_n + \frac{1}{k^2}(-A)^{-3}[-4I + kA + \mathrm{e}^{-kA}(4I + 3kA + k^2 A^2)] (F(u_n, t_n) \nonumber\\
            && + \frac{2}{k^2}(-A)^{-3}[2I - kA - \mathrm{e}^{-kA}(2I + kA)])(F(a_n, t_n + k/2) + F(b_n, t_n + h/2) )  \\
            && + \frac{1}{k^2}(-A)^{-3}[-4I + 3kA - k^2A^2  + \mathrm{e}^{-kA}(4I + kA)] F(c_n, t_n + k) \nonumber
\end{eqnarray}
where, 
\begin{eqnarray*}
    a_n &=& \mathrm{e}^{-kA/2}U_n - A^{-1} \left( \mathrm{e}^{-kA/2} - I \right) F(u_n, t_n) \\
    b_n &=& \mathrm{e}^{-kA/2}U_n - A^{-1} \left( \mathrm{e}^{-kA/2} - I \right)  F(a_n, t_n + k/2)  \\
    c_n &=& \mathrm{e}^{-kA/2}a_n - A^{-1} \left( \mathrm{e}^{-kA/2} - I \right) [2 F(b_n, t_n + k/2) - F(u_n, t_n)] \\
\end{eqnarray*}

\noindent The presence of higher powers of matrix inverses and the matrix exponentials, coupled with the reciprocal of the time step ($k$) in the scheme, makes it computationally expensive. Some authors \cite{asante2025fourth,yousuf2012numerical, yousuf2009efficient}, have derived fourth-order ETDRK4 schemes using fourth-order Pad\'e (0,4) and Pad\'e (2,2) rational functions to approximate the matrix exponential. While Pad\'e (0,4) rational functions are L-acceptable they have complex poles, making their algebraic computations less efficient. On the other hand Pad\'e (2,2) rational functions are only A-acceptable, and thus are unable to efficiently damp spurious oscillations for problems with mismatched initial and boundary conditions without the application of presmoothing steps \cite{asante2025fourth}. To curb these issues, we require a rational function that is both L-acceptable and has real and distinct poles to approximate the matrix exponential.\\   
\noindent
To make the scheme easy to develop, we make the following substitutions to Eq.~\eqref{eqn2.21}, 
\begin{eqnarray}\label{eqn2.22}
    a_n&=& \mathrm{e}^{-\frac{k}{2}A}U_n +\tilde{P}(kA)F(U_n,t_n) \nonumber \\
    b_n&=& \mathrm{e}^{-\frac{k}{2}A}U_n +\tilde{P}(kA)F(a_n,t_n+\frac{k}{2})  \\
    c_n&=& \mathrm{e}^{-\frac{k}{2}A}a_n +\tilde{P}(kA)[2F(b_n,t_n+\frac{k}{2})-F(U_n,t_n)] \nonumber\\
    U_{n+1} &=& \mathrm{e}^{-kA}U_n +P_1(kA) F(U_n,t_n) +2P_2(kA)(F(a_n,t_n+\frac{k}{2})+F(b_n,t_n+\frac{k}{2}))\nonumber \\
    && + P_3(kA)F(c_n,t_n+k)\,, \nonumber
\end{eqnarray}

 with 
\begin{eqnarray}
    P_1(kA) &=&  \frac{1}{k^2}(-A)^{-3}[-4+kA+\mathrm{e}^{-kA}(4+3kA+k^2A^2)]. \\
    P_2(kA) &=&  \frac{1}{k^2}(-A)^{-3}[2-kA-\mathrm{e}^{-kA}(2+kA)].\\
    P_3(kA) &=& \frac{1}{k^2}(-A)^{-3}[-4+3kA-k^2A^2+\mathrm{e}^{-kA}(4+kA)]. \\
    \tilde{P}(kA) &=& - A^{-1}(\mathrm{e}^{-\frac{k}{2}A}-I). \label{ptilde}
\end{eqnarray}

\subsection{Real and distinct pole rational function of the matrix exponential}
 To approximate the matrix exponentials in Eq.~\eqref{eqn2.22},  we consider the fourth-order RDP rational approximation to $\mathrm{e}^z$ developed by  D. A. Voss and A. Q. M. Khaliq in \cite{voss1996time} in 1996. It is of the form 
\begin{equation}\label{eqn2.31}
    R(z) = \frac{N(z)}{D(z)} = \frac{1 + a_1z + a_2z^2 + a_3z^3}{(1-b_1z)(1-b_2z)(1-b_3z)(1-b_4z)}
\end{equation}
where $a_i, b_i \in \mathbb{R}.$ We summarize the corresponding partial fraction decomposition as 
\begin{equation}\label{eqn2.32}
R(z) = \sum_{i=1}^4 \frac{w_i}{1 - b_iz} 
\end{equation}
where the expansion coefficients are estimated in the usual way, $\displaystyle w_i = \lim_{z \longrightarrow \frac{1}{b_k}} (1-b_k z)R(z).$ 

\begin{figure}[ht]
\centering
\includegraphics[width=10cm]{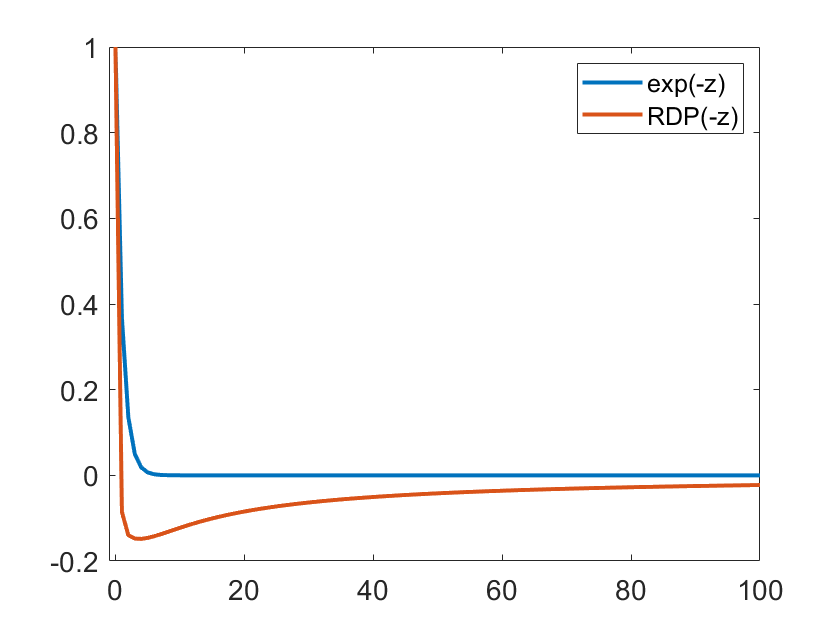}
\caption{Approximation of $\mathrm{e}^{-z}$ with $R(-z)$} 
\label{eqn:figRDP}
\end{figure}

We require the rational approximation to $\mathrm{e}^z$ to be A-acceptable (ie. $|R(z)| < 1$ when $Re(z) < 0, z \in \mathbb{C}$) and L-acceptable (ie, in addition to being A-acceptable, $|R(z)| \longrightarrow 0$ as $z \longrightarrow -\infty$). We also require each $b_i$ to be distinct to allow for parallelization in a multiprocessor environment. By the Maximum modulus theorem, $R(z)$ will be an L-acceptable approximation to $e^z$ if $b_i >0$, and the polynomial $E(y):= |D(\mathrm{i}y)|^2 - |N(\mathrm{i}y)|^2 \geq 0, \forall y \in \mathbb{R}$ \cite{voss1996time}. D. A. Voss and A. Q. M. Khaliq in \cite{voss1996time} developed parameters for the constants in the partial fraction decomposition form of $R(z)$ in Table [\ref{tab_etd4_para}] with the constraint of L-acceptability while maintaining a near optimal error constant. The graph of $R(z)$ using the parameters in Table [\ref{tab_etd4_para}] is shown in Fig.~\ref{eqn:figRDP}.

To derive the ETDRK4RDP scheme, we define a corresponding rational approximation of the matrix exponential. Recall from Cauchy's integral formula for matrix functions that if $f(z)$ is a complex analytic function inside and on a simple closed contour $\Gamma$, and $A$ is a square matrix whose spectrum lies entirely inside $\Gamma$, then the function of the matrix $f(A)$ defined by:
\begin{equation}\label{cauchy1}
    f(A) = \frac{1}{2\pi \mathrm{i}} \int_{\Gamma} f(z)(zI - A)^{-1} \  \mathrm{d} z
\end{equation}
is well-defined by integrating around the eigenvalues of $A$ \cite{gauthier2014lectures, huang2021derivative}. The corresponding matrix notation of the RDP for the matrix exponential is given as 
\begin{equation}\label{eqn2.33}
R(-kA) = (I + a_1 kA + a_2 k^2A^2 + a_3 k^3 A^3)((I + b_1 k A)(I + b_2 k A)(I + b_3 k A)(I + b_4 k A))^{-1}    
\end{equation}
with corresponding partial fraction decomposition of the form 
\begin{equation}\label{eqn2.3311}
  R(-kA) = \sum_{i=1}^4 (w_i)(I + b_i k A)^{-1}. 
\end{equation}

\noindent

\begin{table}
\begin{center}
\begin{tabular}{ c c c c} 
 \hline
  k & $b_k$ & $w_k$& $a_k$ \\
  \hline
  \hline
 1 & 0.4751834017787114 & 20.10707940496431 & \\ 
 2 & 1.0000000000000000 & 0.5229558818011362 & 1.579627631895178\\ 
 3 & 0.3888888888888889 & -15.21083750434353 & 0.303085087930133\\
 4 & 0.7155553412275962 & -4.419197782421921 & -0.324250474367700\\
 \hline
\end{tabular}
\caption{Parameter values: Partial fraction decomposition of the RDP rational function}\label{tab_etd4_para}
\end{center}
\end{table}

\noindent
The RDP approximation to $\mathrm{e}^{-\frac{k}{2}A}$ in matrix notation is given by 
\begin{equation}\label{eqn2.331122}
    R(-kA/2) = \left(I + \frac{k}{2} a_1A + \frac{k^2}{4}a_2A^2 + \frac{k^3}{8}a_3A^3\right)\left((I + \frac{b_1}{2} k A)(I + b_2 k A)(I +\frac{b_3}{2} k A)(I + \frac{b_4}{2} k A)\right)^{-1}
\end{equation}  
\noindent
with partial fraction decomposition of the form 
\begin{equation}\label{eqn2.3331}
    R (-kA/2) = \sum_{i = 1}^{4} w_i\left( I + \frac{b_i}{2}kA \right)^{-1}.
\end{equation}

\subsection{Derivation of the ETDRK4RDP scheme}
We begin the derivation by considering each of the approximations $a_n, b_n, c_n$ and $U_{n+1}$ separately and simplify their results when the matrix exponentials are substituted with the corresponding RDP approximation. Consider $\tilde{P}(kA)$ as defined in Eq.~\eqref{ptilde} and substitute the RDP approximation $R(-kA/2)$ in Eq.~\eqref{eqn2.331122} with $\mathrm{e}^{-kA/2}$. We have 
$$  \tilde{P}(kA) = - A^{-1}(\mathrm{e}^{-\frac{k}{2}A}-I). $$
That is, 
\begin{eqnarray*}
\tilde{P}(kA) &=& - A^{-1}(R(-kA/2)- I)\\
\tilde{P}(kA)  &=& -A^{-1}\left( (I + a_1 kA + a_2 k^2A^2 + a_3 k^3 A^3)\left( \prod_{i=1}^4 (I + \frac{b_i}{2} k A) \right)^{-1} - I \right). 
\end{eqnarray*}
\noindent
Then, we factor out the matrix inverses and simplify algebraically and write the partial fraction decomposition of $\tilde{P}(kA)$ as 

\begin{equation}\label{eqn2.44421}
     \tilde{P}(kA) = \sum_{i = 1}^{4} q_i\left( I + \frac{b_i}{2}kA \right)^{-1},
\end{equation}

where 
$$\displaystyle q_i = \frac{k \left( \mu_1 - 2\frac{\mu_2}{b_i} + 4\frac{\mu_3}{b_i^2} - 8\frac{\mu_4}{b_i^3} \right)}{\prod_{j=1, j\ne i}^{4} \left( 1 - \frac{b_j}{b_i}\right) } \quad \quad i = 1,2,3,4. $$

and, 
\begin{eqnarray*}
    \mu_1 &=& \frac{1}{2} (a_1 - (b_1 +b_2+b_3+b_4)), \quad \quad \mu_2 = \frac{1}{4} (a_2 - (b_1b_2 + b_3b_4 + (b_1 + b_2)(b_3+b_4)) \\
    \mu_3 &=& \frac{1}{8} (a_3 - ((b_1 + b_2)b_3b_4 + (b_3 + b_4)b_1b_2)), \quad \quad \mu_4 = -\frac{1}{16}(b_1b_2b_3b_4). 
\end{eqnarray*}

\noindent
To simplify $a_n$, we substitute the partial fraction decomposition of $\tilde{P}(-kA)$ in Eq.~\eqref{eqn2.44421} into Eq.~\eqref{eqn2.22} and the RDP approximation $R(-kA/2)$ in Eq.~\eqref{eqn2.331122} with $\mathrm{e}^{-kA/2}$. Thus, 
\begin{eqnarray*}
    a_n &=&  \mathrm{e}^{(-kA/2)}U_n +\tilde{P}(kA)F(U_n,t_n)\\
    a_n &=&  R(-kA/2)U_n +\tilde{P}(kA)F(U_n,t_n)\\
    a_n &=&  \sum_{i = 1}^{4} w_i\left( I + \frac{b_i}{2}kA \right)^{-1}U_n +  \sum_{i = 1}^{4} q_i\left( I + \frac{b_i}{2}kA \right)^{-1} F(U_n,t_n).
\end{eqnarray*}

Next, we factor out the sum of the matrix inverse and simplify the resulting term in partial fraction form as 

\begin{equation}\label{eqn2.42}
    a_n = \sum_{i = 1}^4 J_{ia}   
\end{equation}
where, 
$$   J_{ia} =  \left(  I + \frac{b_i}{2}kA \right)^{-1} (w_iu_n - q_iF(u_n, t_n)). $$

\noindent
Similarly for $b_n$,  we substitute the partial fraction decomposition of $\tilde{P}(-kA)$ in Eq.~\eqref{eqn2.44421} into Eq.~\eqref{eqn2.22} and the RDP approximation $R(-kA/2)$ in Eq.~\eqref{eqn2.331122} with $\mathrm{e}^{-kA/2}$.

\begin{eqnarray*}
    b_n &=& \mathrm{e}^{-\frac{k}{2}A}U_n +\tilde{P}(kA)F(a_n,t_n+\frac{k}{2})  \\
    b_n&=& R(-kA/2)U_n +\tilde{P}(kA)F(a_n,t_n+\frac{k}{2})  \\
    b_n &=& \sum_{i = 1}^{4} w_i\left( I + \frac{b_i}{2}kA \right)^{-1}U_n +  \sum_{i = 1}^{4} q_i\left( I + \frac{b_i}{2}kA \right)^{-1} F\left( a_n,t_n+\frac{k}{2} \right).
\end{eqnarray*}

Again, we factor out the matrix inverses and simplify the resulting term in partial fraction form as, 
\begin{eqnarray}\label{eqn2.43}
 b_n = \sum_{i = 1}^4 J_{ib},
\end{eqnarray} 
\noindent
where, 
$$ J_{ib} = \left(  I + \frac{b_i}{2}kA \right)^{-1} (w_iu_n - q_iF(a_n, t_n + k/2)). $$
\noindent

Similarly, for $c_n$ we substitute $\tilde{P}(-kA)$ in Eq.~\eqref{eqn2.44421} into Eq.~\eqref{eqn2.22} and the RDP approximation $R(-kA/2)$ in Eq.~\eqref{eqn2.331122} with $\mathrm{e}^{-kA/2}$. The resulting equation is,

\begin{eqnarray*}
     c_n &=& \mathrm{e}^{-\frac{k}{2}A}a_n +\tilde{P}(kA)[2F(b_n,t_n+ k/2)-F(U_n,t_n)] \\
      c_n&=& R(-kA/2)a_n +\tilde{P}(kA)[2F(b_n,t_n+ k/2)-F(U_n,t_n)] \\
      c_n &=&  \sum_{i = 1}^{4} w_i\left( I + \frac{b_i}{2}kA \right)^{-1}a_n +   \sum_{i = 1}^{4} q_i\left( I + \frac{b_i}{2}kA \right)^{-1} [2F(b_n,t_n+ k/2)-F(U_n,t_n)]. 
\end{eqnarray*}

We again factor out the matrix inverses and simplify the resulting term in partial fraction form as,
\begin{eqnarray}\label{eqn2.44}
 c_n = \sum_{i = 1}^4 J_{ic},
\end{eqnarray} 
where
$$ J_{ic} = \left( I + \frac{b_i}{2}kA \right)^{-1} (w_ia_n - q_i[2F(b_n, t_n + k/2) - F(u_n, t_n)]). $$

\noindent
Finally, to evaluate $U_{n+1}$, we begin by substituting the RDP approximation $R(-kA)$ in Eq.~\eqref{eqn2.3331} with $\mathrm{e}^{-kA}$ and simplify the coefficients of each of the terms in the expression independently. To begin, observe that 
$$ (I + b_1kA)(I + b_2kA)(I + b_3kA)(I + b_4kA) = I + \alpha kA + \beta k^2A^2 + \gamma k^3A^3 + \rho k^4A^4 $$
where 
\begin{eqnarray*}
    \alpha &=& b_1 + b_2 + b_3 + b_4, \quad \quad  \beta = b_1b_2 + b_3b_4 + (b_1 + b_2)(b_3 + b_4)\\
    \gamma &=& (b_1 + b_2)b_3b_4 + (b_3 + b_4)b_1b_2, \quad \quad   \rho = b_1b_2b_3b_4.
\end{eqnarray*}
Consider $P_1(kA)$ and substitute $\mathrm{e}^{-kA}$ with $R(-kA)$ as defined in Eq.~\eqref{eqn2.3331}, thus
\begin{eqnarray*}
    P_1(kA) &=&  \frac{1}{k^2}(-A)^{-3}[-4+kA+ R(-kA)(4+3kA+k^2A^2)]\\
    P_1(kA) &=&  \frac{1}{k^2}(-A)^{-3}[-4+kA+ (I + a_1 kA + a_2 k^2A^2 + a_3 k^3 A^3)\left( \prod_{i=1}^4 (I + b_i k A) \right)^{-1}(4+3kA+k^2A^2)].
\end{eqnarray*}

\noindent
Next, factor out the product of the matrix inverse and simplify the resulting algebraic expression,
$$  P_1(kA)  = (w_1 kI + w_2 k^2A + w_3 k^3A^2)\left(\prod_{i=1}^4(I + b_ikA) \right)^{-1}.$$
Next, we write this rational function in partial fraction form as 
\begin{equation}\label{eqn2.45}
     P_1(kA)= \sum_{i=1}^4  r_i \left(I + b_ikA \right)^{-1},
\end{equation}
where,
$$ r_i = \frac{\left( w_1 - \frac{w_2}{b_i} + \frac{w_3}{b_i^2} \right)k}{\prod_{j=1, j \ne i}^4   \left( 1 - \frac{b_j}{b_i} \right)}, \quad i = 1,2,3,4,$$
with 
$$ w_1 =4\gamma - \beta -a_1 - 3a_2 - 4a_3, \quad w_2 = 4\rho - \gamma -a_2 - 3a_3, \quad w_3 = -\rho - a_3  $$

Similar to how we obtained $P_1(kA)$, we simplify $P_2(kA)$ as, 
\begin{equation}\label{eqn2.46}
    P_2(kA) =  \sum_{i = 1}^4 g_i \left(I + b_i kA \right)^{-1}
\end{equation}  
where, 
$$ g_i = \frac{\left( l_1 - \frac{l_2}{b_i} + \frac{l_3}{b_i^2} \right)k}{\prod_{j=1, j \ne i}^4 \left( 1 - \frac{b_j}{b_i} \right)}, \quad \quad i= 1,2,3,4$$

with, 
$$ l_1 =-(2\gamma - \beta - (2a_3 + a_2)), \quad l_2 = -(2\rho - \gamma - a_3), \quad l_3 = \rho.  $$

Finally, we simplify $P_3(kA)$ in a similar procedure as,
\begin{equation}\label{eqn2.47}
    P_3(kA) =  \sum_{i = 1}^4  h_i \left(I + b_ikA \right)^{-1} 
\end{equation}
where, 
$$ h_i = \frac{\left( m_1 - \frac{m_2}{b_i} + \frac{m_3}{b_i^2} - \frac{m_4}{b_i^3}\right)k}{\prod_{j=1, j \ne i}^4 \left( 1 - \frac{b_j}{b_i} \right)}, \quad \quad i= 1,2,3,4 $$

with, 
$$m_1 = -(-4\gamma +3\beta -\alpha+(a_2 + 4a_3)), \quad m_2 = -(-4\rho +3\gamma -\beta + a_3), \quad m_3 = -(3\rho - \gamma), \quad m_4 = \rho.  $$

Substituting the simplified coefficients $P_1(kA) , P_2(kA)$ and $P_3(kA)$ in Eq.~\eqref{eqn2.45}, Eq.~\eqref{eqn2.46}, and Eq.~\eqref{eqn2.47} respectively, into the final approximation $U_{n+1}$ in Eq.~\eqref{eqn2.21}, we have that

\begin{eqnarray*}
     U_{n+1} &=& \left(\sum_{i=1}^4 (w_i)(I + b_i k A)^{-1}\right) U_n + \left(\sum_{i=1}^4  r_i \left(I + b_ikA \right)^{-1} \right) F(U_n,t_n) \nonumber \\
     && +2  \left(\sum_{i=1}^4  g_i \left(I + b_ikA \right)^{-1} \right) (F(a_n,t_n+ k/2)+F(b_n,t_n+ k/2) +  \left(\sum_{i=1}^4  h_i \left(I + b_ikA \right)^{-1} \right) F(c_n,t_n+k). \nonumber
\end{eqnarray*}

\noindent
We then group the terms with similar quotients and summarize as follows 
\begin{equation*}
    U_{n+1} = \sum_{i = 1}^4 [ \left(I + b_ikA \right)^{-1} ( w_iu_n + r_iF(u_n, t_n) + 2g_i[F(a_n, t_n + k/2) + F(b_n, t_n + k/2)] + h_iF(c_n, t_n + k) ) ].
\end{equation*}
\noindent
Rewriting the solution $U_{n+1}$ as a linear system, we solve for $Y_i$ for $i = 1,2,3,4$ in the following systems 
\begin{equation*}
    (I + b_ikA) Y_i = w_iu_n + r_iF(u_n, t_n) + 2g_i[F(a_n, t_n + k/2) + F(b_n, t_n + k/2)] + h_iF(c_n, t_n + k)   
\end{equation*}
\noindent
and obtain $U_{n+1}$ as 
$$ U_{n+1} = \sum_{i=1}^4 Y_i$$

\subsection{Implementation of the ETDRK4RDP algorithm}
Putting all the results together, the ETDRK4RDP scheme can be implemented in serial with the following steps:
\begin{enumerate}
    \item Solve the system $\displaystyle \left( I + \frac{b_i}{2}kA \right)J_{ia} = (w_iU_n - q_iF(U_n, t_n))$ for $i = 1,2,3,4$.\label{A1}
    \item Solve for $a_n$ as $\displaystyle a_n = \sum_{i=1}^4 J_{ia}$.
    \item Solve the system $\displaystyle \left( I + \frac{b_i}{2}kA \right)J_{ib} = (w_iU_n - q_jF(a_n, t_n +k/2))$ for $i = 1,2,3,4$.\label{B1}
    \item Solve for $b_n$ as $\displaystyle b_n = \sum_{i=1}^4 J_{ib}$.
    \item Solve the system $\displaystyle \left( I + \frac{b_i}{2}kA \right)J_{ic} = (w_ia_n - q_j[2F(b_n, t_n + k/2) - F(U_n, t_n)])$ for $i = 1,2,3,4$.\label{C1}
    \item Solve for $c_n$ as $\displaystyle c_n = \sum_{i=1}^4 J_{ic} $.
    \item Solve the system 
    $$ (I + b_jkA) Y_j = w_jU_n + r_jF(U_n, t_n) + 2g_j[F(a_n, t_n + k/2) + F(b_n, t_n + k/2)] + h_jF(c_n, t_n + k) $$
    \noindent
    for $j = 1,2,3,4$. \label{D1}
    \item Solve for $U_{n+1}$ as $U_{n+1} = Y_1 + Y_2 + Y_3 + Y_4$.
\end{enumerate}

\noindent
To implement the algorithm in a multiprocessor environment, we use the UMFPACK package (see \cite{Davis2025suitesparse}) to solve the sparse linear system of equations with the matrices stored in compressed column storage format (see \cite{Hanyk2014mUMFPACK}). A similar implementation has been used in \cite{asante2020second} to achieve similar results. The algorithm is implemented in a Fortran interface, using OpenMP with 4 processors. At each time step, 4 independent processors are tasked to solve each of the systems in step (\ref{A1}) for $J_{ia}, i = 1,2,3,4$ and the results compiled to find $a_n$. Similar implementations are performed in steps (\ref{B1}), (\ref{C1}), and (\ref{D1}) to obtain $b_n, c_n$ and the solution $U_{n+1}$ respectively.

\section{Numerical Experiments}
\label{section4}
To investigate the order of accuracy and efficiency of the proposed scheme, we tested the performance of the scheme on a variety of test problems with various initial and boundary conditions. We compare the accuracy of the new scheme with existing fourth-order ETDRK4P22 and ETDRK4P04 schemes developed by \cite{yousuf2009efficient} and \cite{asante2025fourth}, respectively. The ETDRK4P22 scheme developed by M. Yousuf for semilinear parabolic problems uses the Pad\'e(2,2) rational function, which is not L-acceptable while the ETDRK4P04 uses the Pad\'e(0, 4) scheme to approximate the matrix exponentials in the ETDRK4 scheme developed by Cox and Matthews \cite{mathews2004numerical}. We evaluate the performance of the schemes on problems in two spatial dimensions as well as for scalar and systems of reaction–diffusion equations with homogeneous Dirichlet and Neumann boundary conditions. The error, $E(k)$ at time step $k$, is measured with the $L_{\infty}-$ norm and the estimated order of convergence is calculated using the formula     
$$ p = \frac{\log[E(k)/E(k/2)]}{\log 2}. $$
\noindent
For problems with known exact solution, we estimate the error as $E(k) = \|U(k) - \hat{U}(k) \|_{\infty}$, where $\hat{U}(k)$ and $U(k)$ denote the approximate grid solution and exact solutions, respectively at the final time $T$ using a temporal step size of $k$. The spatial derivatives are approximated with a fourth-order finite difference scheme on a mesh with spacing $h$ described in Section \ref{section2}. By setting $k=h$, we were able to derive the correct order of convergence of the fully discretized system. For problems with unknown exact solution, we use the numerical solution on the next level of refinement as the reference solution for computing errors and the order of convergence, $\tilde{p}$. We only evaluated the convergence of the time discretization by fixing the spatial mesh size, $h$. The errors are computed by using the numerical solution on the next fine mesh as the reference solution. For a grid refinement study with temporal step sizes $k, k/2, k/4, k/8 $ we use the formula
$$ \tilde{p} = \frac{\log[\tilde{E}(k)/\tilde{E}(k/2)]}{\log 2}  $$
where, $\tilde{E}(k) = \|\hat{U}(k) - \hat{U}(k/2) \|$. It's been shown in \cite{leveque2007finite} that the estimation of the error this way determines the correct order of convergence of the scheme for sufficiently small $k$. Estimating the error and testing the order of accuracy by this approach only confirms that the code is converging to some function with the desired rate \cite{leveque2007finite}. The code may be converging very nicely, but to the wrong function. So, test problems with known exact solution are critical for ensuring that the scheme is converging to the correct solution.

Numerical experiments were run in MATLAB R2023b on a Dell Latitude 5440 with 13th Gen Intel(R) Core(TM) i7-1355U 1.70 GHz with 32 GB RAM. To validate the assertion that the proposed scheme has a better efficiency when parallelized, we run the algorithm in Fortran (using OpenMP with 4 processors). The timings are reported accordingly, where we used a PC with 32 $\times$ 12th Gen Intel Core i9-14900 with 32 GiB of RAM and 1.0TB of memory with Ubuntu 22.04 as the operating system. Numerical simulations are performed on four different test problems with homogeneous Dirichlet and Neumann boundary conditions. Two of the problems are scalar with known exact solutions, a nonlinear problem with mismatched initial and boundary data, and finally a system of nonlinear reaction–diffusion equations, called the Brusselator model \cite{zegeling2004adaptive}. 


\subsection{Model problem with a Dirichlet boundary condition}\label{Dirichlet}
We investigate the performance of the proposed scheme by considering the following two-dimensional reaction–diffusion equation with corresponding homogeneous Dirichlet boundary and initial condition
\begin{equation}
    \begin{cases}
        \displaystyle \frac{\partial u}{\partial t}  = \Delta u - u, \quad \frac{-\pi}{2} < x, y < \frac{\pi}{2}, \quad t \in [0, T] \\
        u(x, y, 0) = \cos(x) \cos(y)
    \end{cases}
\end{equation}
\noindent
The exact solution is given by $u(x, y, t) = \mathrm{e}^{-3t}\cos(x) \cos(y). $ We discretize each dimension of the spatial domain with $m+2$ nodes with spatial step size $\displaystyle h = \frac{\pi}{m +1}$. The spatial discretization is performed as discussed in Section \ref{section2}. The proposed ETDRK4RDP scheme converges with fourth-order accuracy (See Fig.~\ref{Fig. 3.11}A and Table~\ref{table 1.225}). For all time steps examined, the serial implementation of ETDRK4RDP scheme is less accurate compared to  ETDRK4P04 scheme and requires more CPU time (See Fig.~\ref{Fig. 3.11} and Table~\ref{table 1.225}). It should be mentioned that the ETDRK4RDP scheme performs twice the number of linear solves as the ETDRK4P04 scheme but only increases the CPU time by 26\%, suggesting a great potential to improve the CPU time through parallelization. 

\begin{figure}[ht]
\centering
\includegraphics[width=12cm]{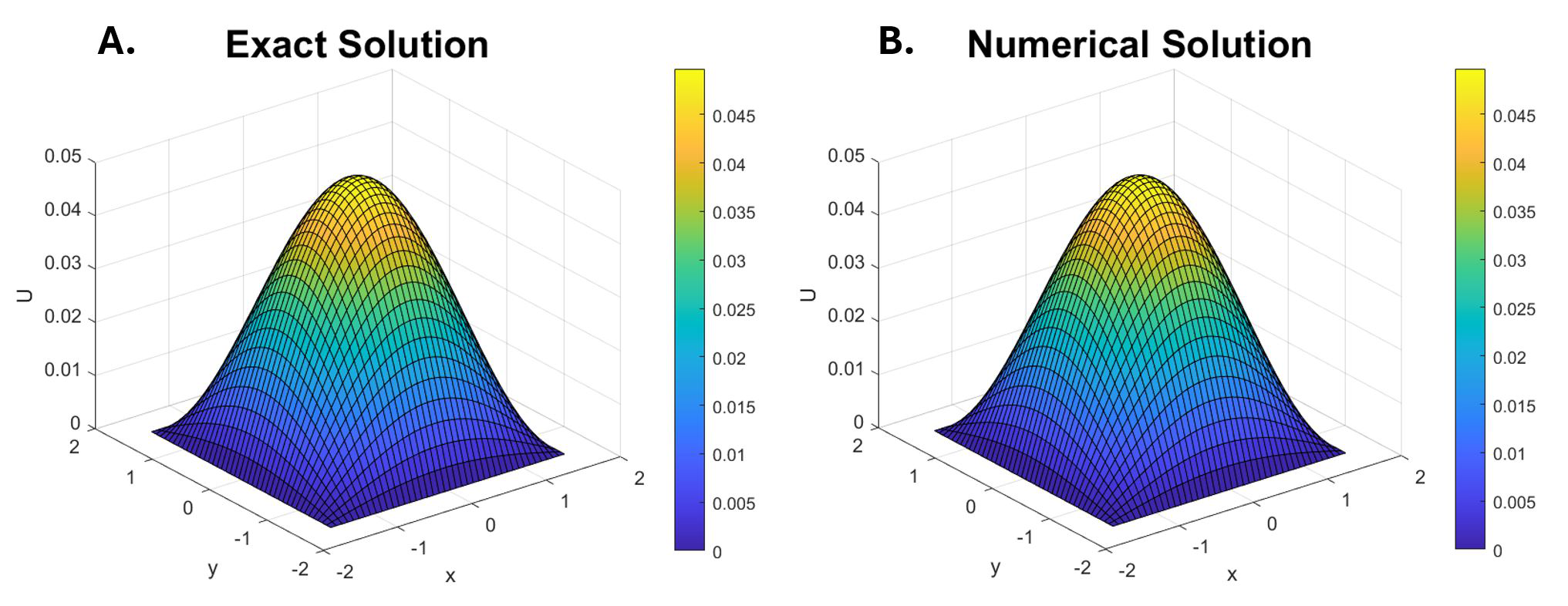}
\caption{Solution profile for (A) Exact solution and (B) Numerical solution using ETDRK4RDP for test-problem \ref{Dirichlet}} 
\label{Fig. 3.11112}
\end{figure}

As seen in Table~\ref{table 1.2232}, our serial implementation of ETDRK4RDP in Fortran using UMFPACK is about three times faster than the MATLAB implementation. We obtain an additional speed-up by a factor of two using the parallelized Fortran program with four processors. Comparing our parallelized ETDRK4RDP in Fortran with the serial ETDRK4P04 in MATLAB, we achieve speed-up by a factor of 5 at the finest time step $k=0.0125)$. 

\begin{figure}[ht]
\centering
\includegraphics[width=12cm]{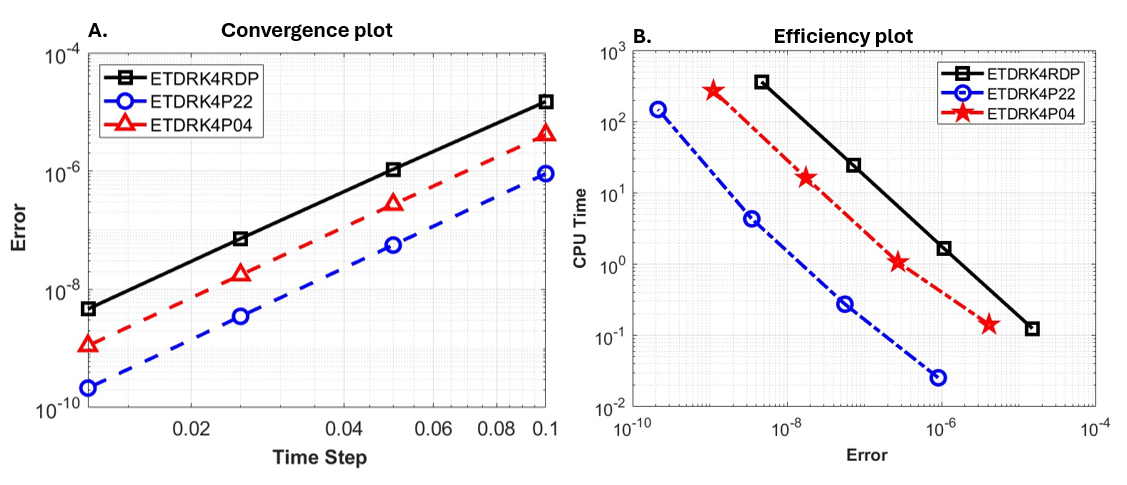}
\caption{Convergence and efficiency plots for Example \ref{Dirichlet}. (A) Convergence and (B) Efficiency of the new fourth-order ETDRK4RDP scheme compared with
the existing fourth-order ETDRK4P22, and ETDRK4P04 schemes.} 
\label{Fig. 3.11}
\end{figure}

\begin{table}[ht]
 \begin{center}
 \resizebox{\textwidth}{!}{
\begin{tabular}{ p{1cm}|p{1cm}||p{2cm}|p{1cm}|p{2cm}||p{2cm}|p{1cm}|p{2cm} }
 \hline
 \multicolumn{8}{c}{\quad ETDRK4RDP \quad \quad \quad \quad \quad \quad \quad \quad \quad ETDRK4P04} \\
 \hline
 \hline
   $k$      & $h$     & Error                 & $p$    & Cpu Time (sec) & Error           & $p$    & Cpu Time (sec) \\
   \hline
   0.10   & 0.08  & $1.50 \times 10^{-5}$ & $-$ & 0.0659   & $4.12 \times 10^{-6}$ & $-$ & 0.0692 \\       
   0.05   & 0.04  & $1.07 \times 10^{-6}$ & 3.80 & 0.8084   & $2.73 \times 10^{-7}$ & 3.92 & 0.7241 \\
   0.025  & 0.02  & $7.23 \times 10^{-8}$ & 3.89 & 12.0560  & $1.76 \times 10^{-8}$ & 3.96 & 9.7913 \\   
   0.0125 & 0.010 & $4.66 \times 10^{-9}$ & 3.96 & 178.1067 & $1.11 \times 10^{-9}$ & 3.98 & 140.2890 \\
 \hline
\end{tabular}
  }
  \caption{\label{table 1.225} Convergence table showing fourth-order convergence of the new ETDRK4RDP scheme compared with the existing ETDRK4P04 scheme for the model test problem in Example \ref{Dirichlet} with $T =1$.}
\end{center}
\end{table}

\begin{table}[!ht]
\begin{center}
\begin{tabular}{ccrrr}
\hline
$k$ & $h$ & Matlab & Fortran (1 processor) & Fortran (4 processors)\\
\hline
\hline
$0.1$ & $0.08$ & 0.066 & 0.03& 0.01\\
$0.05$ & $0.04$  & 0.808& 0.38& 0.14\\
$0.025$ & $0.02$ & 12.056& 4.84& 2.39\\
$0.0125$ & $0.01$ & 178.107& 53.83& 27.50 \\
\hline
\end{tabular}
\caption{\label{table 1.2232}Comparison of CPU times in seconds for the ETD4RK4RDP method implemented in Matlab and Fortran (serial and parallelized version).}
\end{center}
\end{table}

\subsection{Model problem with Neumann boundary}\label{Neumann}
\begin{figure}[ht]
\centering
\includegraphics[width=12cm]{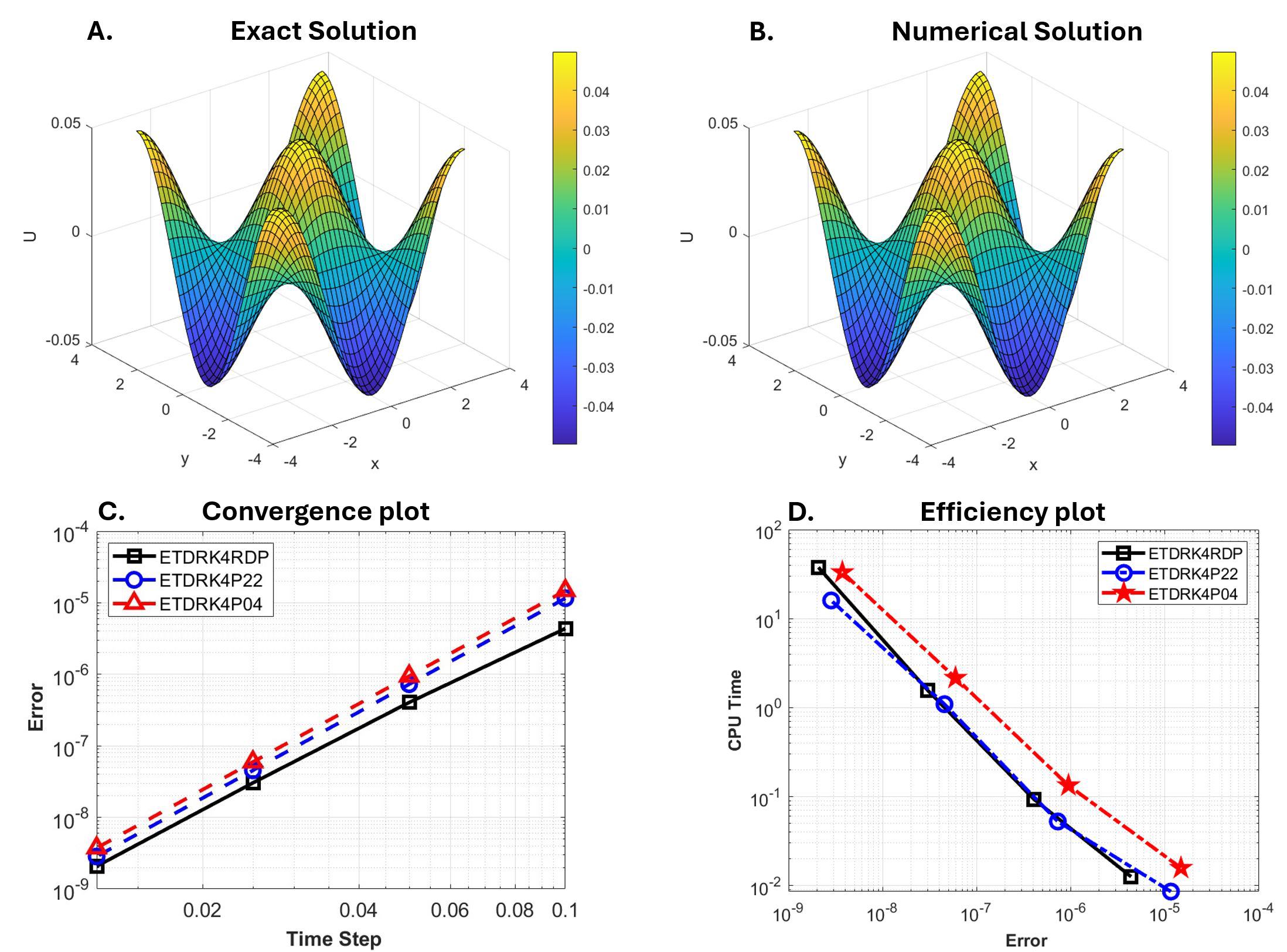}
\caption{Convergence and efficiency plots for Example \ref{Neumann}. (A) Convergence and (B) Efficiency of the new fourth-order ETDRK4RDP scheme compared with the existing fourth-order ETDRK4P22, and ETDRK4P04 schemes.} 
\label{Fig 3.12}
\end{figure}

Next, we consider the following scalar reaction-diffusion system with homogeneous Neumann boundary conditions given by
\begin{equation}
    \begin{cases}
        \displaystyle \frac{\partial u}{\partial t}  = \Delta u - u, \quad -\pi < x, y < \pi, \quad t \in [0, T] \\
        u(x, y, 0) = \cos(x) \cos(y)
    \end{cases}
\end{equation}
\noindent
The corresponding exact solution is
$u(x, y, t) = \mathrm{e}^{-3t}\cos(x) \cos(y)$. We discretize each dimension of the spatial domain with $m+2$ nodes and set a spatial step size $\displaystyle h = \frac{2\pi}{m+1}$. The spatial discretization is done as described in Section \ref{section2}. Fig.~\ref{Fig 3.12}A and Fig.~\ref{Fig 3.12}B depict a plot of the exact solution and numerical solution using ETDRK4RDP, respectively. Here we see a display of how close the numerical solution is to the exact one. 

Again, our empirical convergence analysis in Table~\ref{table 3.111} and convergence plot in Fig.~\ref{Fig 3.12}C show that the ETDRK4RDP scheme converges to the correct solution with order four. Whereas the errors in this case are comparable with the competing ETDRK4P04 scheme, our serial implementation of ETDRK4RDP is not as fast. However, our parallelized version of the scheme is about five times faster than the serial ETDRK4P04 scheme (See Table~\ref{table 1.312}).

\begin{table}[ht]
 \begin{center}
  \resizebox{\textwidth}{!}{
\begin{tabular}{ p{1cm}|p{1cm}||p{2cm}|p{1cm}|p{2cm}||p{2cm}|p{1cm}|p{2cm} }
 \hline
 \multicolumn{8}{c}{ETDRK4RDP \quad \quad \quad \quad \quad \quad \quad \quad \quad \quad \quad ETDRK4P04} \\
 \hline
 \hline
   $k$      & $h$    & Error                 &$p$    & Cpu Time (sec) & Error           & $p$    & Cpu Time (sec) \\
   \hline
   0.10   & 0.16 & $4.37 \times 10^{-6}$ & $-$ & 0.0041   & $1.48 \times 10^{-5}$ & $-$ & 0.0102 \\       
   0.05   & 0.08 & $4.05 \times 10^{-7}$ & 3.53 & 0.1032  & $9.43 \times 10^{-7}$ & 3.97 & 0.0846 \\
   0.025  & 0.04 & $3.03 \times 10^{-8}$ & 3.74 & 1.4156  & $5.95 \times 10^{-8}$ & 3.99 & 1.1445 \\   
   0.0125 & 0.02 & $2.07 \times 10^{-9}$ & 3.87 & 23.9256 & $3.74 \times 10^{-9}$ & 3.99 & 17.7786 \\
 \hline
\end{tabular}
 }
\end{center}
\caption{\label{table 3.111} Convergence table showing fourth-order convergence of the new ETDRK4RDP scheme compared with the existing ETDRK4P04 scheme for the model test
problem in Example \ref{Neumann} with $T =1$ }
\end{table}

\begin{table}[!ht]
\begin{center}
\begin{tabular}{ccrrr}
\hline
$k$ & $h$ & Matlab & Fortran (1 processor) & Fortran (4 processors)\\
\hline
\hline
$0.1$ & $0.16$ & 0.0041 & 0.004 & 0.0026\\
$0.05$ & $0.08$  & 0.1032& 0.079 & 0.021\\
$0.025$ & $0.04$ & 1.4156 & 0.927 & 0.232\\
$0.0125$ & $0.02$ & 23.93 & 15.874 & 3.969 \\
\hline
\end{tabular}
\caption{\label{table 1.312}Comparison of CPU times in seconds for the ETDRK4RDP method implemented in Matlab and Fortran (serial and parallelized version).}
\end{center}
\end{table}

\subsection{Nonlinear problem with mismatched initial and boundary data}\label{mismatched}
\begin{figure}[ht]
\centering
\includegraphics[width=12cm]{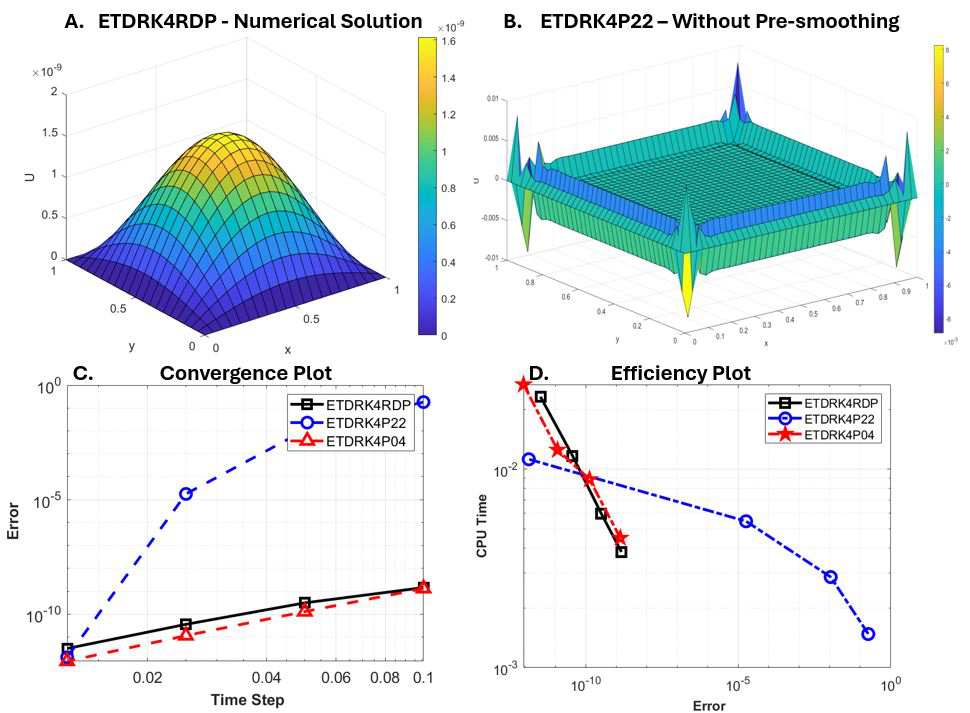}
\caption{Convergence and efficiency plots for Example \ref{mismatched}. (A) Numerical-solution using ETDRK4RDP, (B) ETDRK4P22-solution without Pre-smoothing (C) Convergence Plot (D) Efficiency of the new fourth-order ETDRK4RDP scheme compared with the existing fourth-order ETDRK4P22, and ETDRK4P04 schemes.} 
\label{Fig 3.14}
\end{figure}

Next, we consider our first nonlinear reaction-diffusion equation with homogeneous Dirichlet conditions set at the boundaries. The model problem is the two-dimensional analogue of the one-dimensional reaction–diffusion equation with Michaelis–Menten enzyme kinetics \cite{muller2015methods, bhatt2015locally, cherruault1990stability}. The model is the simplest case of enzyme kinetics, applied to an enzyme-catalysed reaction involving the transformation of one substrate into one product \cite{muller2015methods}. It is described mathematically as 
\begin{equation}
    \begin{cases}
        \displaystyle \frac{\partial u}{\partial t} = d \left(\frac{\partial^2u}{\partial x^2} +\frac{\partial^2u}{\partial y^2} \right) - \frac{u}{(1+u)}\,,\quad 0<x\,,\,y<1\,,\quad t>0 \\
        u(x, y, 0) = 1, \quad 0 \leq x, y \leq 1.
    \end{cases}
\end{equation}
\noindent
where the diffusion coefficient $d = 1$. The problem has no known exact solution. We observe that the initial condition is discontinuous at the boundary where the Dirichlet boundary conditions are applied. Such mismatches of the initial and boundary conditions are known to generate spurious oscillations, which, if not properly damped at the initial stages of the simulation, can significantly affect the accuracy of solutions \cite{asante2025fourth}. Without pre-smoothing, Pad\'e schemes which are not L-acceptable (ie. ETDRK4P22) fail to damp out spurious oscillations (See Fig.~\ref{Fig 3.14}B),  and struggle to maintain fourth-order accuracy Table~\ref{tab etdrk4rdp_conv_eg4_1} and Fig.~\ref{Fig 3.14}D). This phenomenon is also established in \cite{asante2025fourth}. The L-acceptability of the proposed scheme allows the smoothing of the boundary (see Fig.~\ref{Fig 3.14}A), ensuring that the scheme converges with fourth-order accuracy (Table~\ref{table 3.313}). In Fig.~\ref{Fig 3.14}C, we see the inconsistency in the convergence results of ETDRK4P22 scheme since we performed no pre-smoothing. 

Here, the serial version of the proposed scheme is slightly faster than the ETDRK4P04 scheme across all time steps (Table~\ref{table 3.313}). Interestingly, in Table~\ref{table 1.4} we see that the parallelized ETDRK4RDP scheme does not significantly speed up the serial implementation. This could be a result of the short time duration of the simulation. We are still investigating this observation.

\begin{table}[ht]
\begin{center}
 \resizebox{\textwidth}{!}{
\begin{tabular}{ p{1cm}|p{1cm}||p{2cm}|p{1cm}|p{2cm}||p{2cm}|p{1cm}|p{2cm} }
 \hline
 \multicolumn{8}{c}{ETDRK4RDP \quad \quad \quad \quad \quad \quad \quad \quad \quad ETDRK4P04} \\
 \hline
 \hline
   $k$      & $h$    & Error                 & $p$    & Cpu Time (sec) & Error           & $p$    & Cpu Time (sec) \\
   \hline
   0.10   & 0.05 & $1.45 \times 10^{-9}$ & $-$ & 0.0038   & $1.3 \times 10^{-9}$ & $-$ & 0.0045 \\       
   0.05   & 0.05 & $3.1 \times 10^{-10}$ & 2.21 & 0.0060  & $1.3 \times 10^{-10}$ & 3.35 & 0.0089 \\
   0.025  & 0.05 & $3.6 \times 10^{-11}$ & 3.41 & 0.0117  & $1.2 \times 10^{-11}$ & 3.48 & 0.0125 \\   
   0.0125 & 0.05 & $3.2 \times 10^{-12}$ & 3.60 & 0.0232 & $8.9 \times 10^{-13}$ & 3.72 & 0.0267\\
 \hline
\end{tabular}
}
\end{center}
\caption{Convergence table showing fourth-order convergence of the new ETDRK4RDP scheme compared with the existing ETDRK4P04 scheme for the model test problem in Example \ref{mismatched} with $T = 1.$}
\label{table 3.313}
\end{table}

\begin{table}[!ht]
\begin{center}
\begin{tabular}{ccrrr}
\hline
$k$ & $h$ & Matlab & Fortran (1 processor) & Fortran (4 processors)\\
\hline
\hline
$0.1$ & $0.05$ & 0.0038 & 0.0036 & 0.002 \\         
$0.05$ & $0.05$  & 0.0060 & 0.0069 & 0.004 \\
$0.025$ & $0.05$ & 0.0117 & 0.0182 & 0.008\\
$0.0125$ & $0.05$ & 0.0232 & 0.0536 & 0.016 \\
\hline
\end{tabular}
\caption{\label{table 1.4}Comparison of CPU times in seconds for the ETD4RK4RDP method implemented in Matlab and Fortran (serial and parallelized version) for test problem \ref{mismatched}.}
\end{center}
\end{table}

\subsection{The Brusselator model in 2D: A system of reaction-diffusion equations }
\begin{figure}[ht]
\centering
\includegraphics[width=12cm]{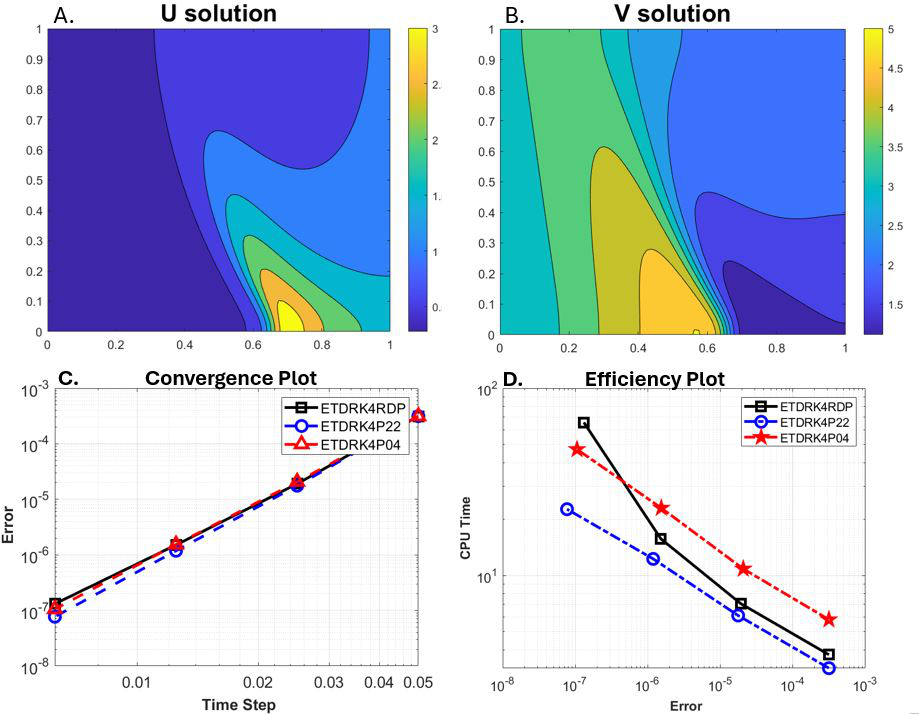}
\caption{Convergence and efficiency plots for Example \ref{Neumann}. (A) U-solution and (B) V-Solution (C) Convergence Plot (D) Efficiency of the new fourth-order ETDRK4RDP scheme compared with the existing fourth-order ETDRK4P22 and ETDRK4P04 schemes.} 
\label{3.15}
\end{figure}

Finally, we investigate the performance of the proposed scheme ETDRK4RDP in simulating a system of non-linear reaction–diffusion equations, with the Brusselator model \cite{zegeling2004adaptive} as a test case. It is a mathematical model of chemical reaction dynamics developed by Ilya Prigogine and colleagues in the 1960s. The model is a nonlinear two-component reaction system that exhibits complex behavior, including oscillations, spatial patterns, and chaos. The nonlinearity in the model coupled with oscillating dynamics, makes numerical computations challenging. The model is given by

\begin{eqnarray*}
u_t &=& \epsilon_1 \Delta u + A + u^2v-(B+1)u\,,\\
v_t &=& \epsilon_2 \Delta v + Bu - u^2v\,,
\end{eqnarray*}
with, $t\in (0,2]$, diffusion coefficients $\epsilon_1=\epsilon_2=2\cdotp 10^{-3}$, and chemical parameters $A=1$, $B=3.4$ in 2D. Homogeneous Neumann boundary conditions are imposed at the boundary of the domain. The initial conditions are 
$$ u(x, y, 0) = \frac{1}{2} + y, \quad v(x, y, 0) = 1 + 5x, \quad \quad 0 < x, y < 1. $$
The problem has no exact solution, hence, we consider investigating the order and accuracy with a fine grid solution. The spatial derivatives are discretized as described in Section \ref{section2}. The numerical solutions, convergence, and efficiency plots are shown in Fig.~\ref{3.15}. 
\noindent
The proposed ETDRK4RDP scheme converges with fourth-order accuracy. ETDRK4RDP performs faster in CPU time than ETDRK4P04 across all time steps and spatial resolutions. Higher efficiencies by almost a factor of 6 are obtained by implementing the ETDRK4RDP scheme in a multiprocessor environment,see Table.~\ref{table 4.4jh}. 

\begin{table}[ht]
\begin{center}
 \resizebox{\textwidth}{!}{
\begin{tabular}{ p{2cm}|p{1cm}||p{2cm}|p{1cm}|p{2cm}||p{2cm}|p{1cm}|p{2cm} }
 \hline
 \multicolumn{8}{c}{ETDRK4RDP \quad \quad \quad \quad \quad \quad \quad \quad \quad \quad \quad ETDRK4P04} \\
 \hline
 \hline
   $k$      & $h$    & Error                 & $p$    & Cpu Time (sec) & Error           & $p$    & Cpu Time (sec) \\
   \hline
   0.05   & 0.0125 & $3.14 \times 10^{-4}$ & $-$ & 3.060   & $3.14 \times 10^{-4}$ & $-$ & 5.81 \\       
   0.025   & 0.0125 & $1.91 \times 10^{-5}$ & 4.04 & 6.080  & $2.07 \times 10^{-5}$ & 3.92 & 10.84 \\
   0.0125  & 0.0125 & $1.51 \times 10^{-6}$ & 3.76 & 12.925  & $1.53 \times 10^{-6}$ & 3.75 & 22.21 \\   
   0.00625 & 0.0125 & $1.31 \times 10^{-7}$ & 3.63 & 32.904 & $1.06 \times 10^{-7}$ & 3.86 & 47.16 \\
 \hline
\end{tabular}
 }
\end{center}
\end{table}


\begin{table}[!ht]
\begin{center}
 \resizebox{0.8 \textwidth}{!}{
\begin{tabular}{ccrrr}
\hline
$k$ & $h$ & Matlab & Fortran (1 core) & Fortran (4 cores)\\
\hline
\hline
$0.05$ & $0.0125$ & 3.060 & 0.8319    &   0.2087 \\   
$0.025$ & $0.0125$  & 6.080 & 1.6687  &  0.4173 \\
$0.0125$ & $0.0125$ & 12.925 & 3.2763  &  0.8197 \\
$0.00625$ & $0.0125$ & 32.904 & 6.5789 &  1.6455\\
\hline
\end{tabular}
}
\caption{\label{table 4.4jh}Comparison of CPU times in serial and multiprocessor environment}
\end{center}
\end{table}

\section{Conclusion and Future Work}\label{section5}

In this article, we have developed a fourth-order exponential time differencing scheme of Runge-Kutta (ETDRK4) type that uses a real and distinct pole (RDP) rational function to approximate the matrix exponential. We have established the order of convergence using two linear problems with exact solutions and two non-linear problems with Dirichlet and Neumann boundary conditions in 2D. One of our non-linear problems is made up of a system with no exact solution. For the non-linear problems, the scheme is more computationally efficient than similar competing schemes when executed with single processors. The ETDRK4RDP scheme performs twice as many linear solves (16 solves) as compared to the competing Pad\'e schemes (8 linear solves), yet compares efficiently. We expect the extra savings to persist for problems with periodic boundary conditions. The superior efficiency of the parallelized ETDRK4RDP scheme is justified after running simulations in Fortran with OpenMP using 4 independent processors. We gain speed-ups by almost a factor of 6 when we implement the new scheme in a multiprocessor environment.

To improve the efficiency of the proposed ETDRK4RDP scheme, we plan to develop a dimensional splitting version that uses only the system matrices $A_1$ and $A_2$ introduced in Section \ref{section2}. In 2020, Asante-Asamani, Wade, and Bruce \cite{asante2016dimensional} reformulated Eq.~\eqref{eqn1.2} to involve only the system matrices $A_1$ and $A_2$ and used this to develop a dimensional splitting version of the ETDRK4 scheme developed by Cox and Mathews in 2002 \cite{mathews2004numerical}. By approximating the matrix exponentials in the reformulated ETDRK4 scheme with their respective RDP rational approximations, we develop the dimensional splitting version called ETDRK4RDP-DS. Second-order versions of this scheme have been shown to gain CPU-time speed-ups by a factor of 20 \cite{asante2025fourth}. We also plan to prove rigorously the fourth-order accuracy of the proposed scheme, extend the test problems to RDEs with advection components, and establish all of its error estimates.

\newpage
\subsection{Appendix A: Supplementary data}
\begin{table}[ht]
\begin{tabular}{ p{1cm}|p{1cm}||p{2cm}|p{1cm}|p{2cm}||p{2cm}|p{1cm}|p{2cm} }
 \hline
 \multicolumn{8}{c}{ETDRK4RDP \quad \quad \quad \quad \quad \quad \quad \quad \quad ETDRK4P22} \\
 \hline
 \hline
   $k$      & $h$     & Error                 & $p$    & Cpu Time (sec) & Error                 & $p$    & Cpu Time (sec)  \\
   \hline
   0.10   & 0.08  &  $1.50 \times 10^{-5}$ & $-$ & 0.0659  & $9.07 \times 10^{-7}$ & $-$ & 0.02503  \\       
   0.05   & 0.04  &  $1.07 \times 10^{-6}$ & 3.80 & 0.8084   & $5.62 \times 10^{-8}$ & 4.01 & 0.30115 \\
   0.025  & 0.02  &  $7.23 \times 10^{-8}$ & 3.89 & 12.0560  & $3.50 \times 10^{-9}$ & 4.01 & 4.65940 \\   
   0.0125 & 0.010 &  $4.66 \times 10^{-9}$ & 3.96 & 178.1067 & $2.1 \times 10^{-10}$ & 4.03 & 69.6025 \\
 \hline
\end{tabular}
\caption{Convergence table showing fourth-order convergence of the new ETDRK4RDP compared with the existing ETDRK4P22 scheme for the model test problem with exact solution at final time $T = 1$.}
\label{tab:etdrk4p22_conv_eg1}
\end{table}

\begin{table}[ht]
\begin{tabular}{ p{1cm}|p{1cm}||p{2cm}|p{1cm}|p{2cm}||p{2cm}|p{1cm}|p{2cm} }
 \hline
 \multicolumn{8}{c}{ETDRK4RDP \quad \quad \quad \quad \quad \quad \quad \quad \quad ETDRK4P22} \\
 \hline
 \hline
  $ k$      & $h$    & Error                 &  $p$    & Cpu Time (sec) & Error           & $p$    & Cpu Time (sec) \\
   \hline
   0.10   & 0.16 & $4.37 \times 10^{-6}$ & $-$ & 0.0041   & $1.16 \times 10^{-5}$ & $-$ & 0.0056 \\       
   0.05   & 0.08 & $4.05 \times 10^{-7}$ & 3.53 & 0.1032  & $7.27 \times 10^{-7}$ & 3.99 & 0.0332 \\
   0.025  & 0.04 & $3.03 \times 10^{-8}$ & 3.74 & 1.4156  & $4.54 \times 10^{-8}$ & 4.00 & 0.5576 \\   
   0.0125 & 0.02 & $2.07 \times 10^{-9}$ & 3.87 & 23.9256 & $2.84 \times 10^{-9}$ & 4.00 & 8.6632 \\
 \hline
\end{tabular}
\caption{Convergence table showing fourth-order convergence of the new ETDRK4RDP compared with the existing ETDRK4P22 and ETDRK4P22 schemes for the model test-problem 2 with exact solution at final time $T = 1$.}
\label{tab:etdrk4rdp_conv_eg2}
\end{table}

\begin{table}[ht]
\begin{tabular}{ p{1cm}|p{1cm}||p{2cm}|p{1cm}|p{2cm}||p{2cm}|p{1cm}|p{2cm} }
 \hline
 \multicolumn{8}{c}{ETDRK4RDP \quad \quad \quad \quad \quad \quad \quad \quad \quad ETDRK4P22} \\
 \hline
 \hline
   $k$      & $h$    & Error                 & $p$    & Cpu Time (sec) & Error           & $p$    & Cpu Time (sec) \\
   \hline
   0.10   & 0.165 & $2.33 \times 10^{-5}$ & $-$ & 0.0037   & $1.93 \times 10^{-6}$ & $-$ & 0.0015 \\       
   0.05   & 0.165 & $2.01 \times 10^{-6}$ & 3.63 & 0.0073  & $1.17 \times 10^{-7}$ & 4.05 & 0.0028 \\
   0.025  & 0.165 & $1.52 \times 10^{-7}$ & 3.73 & 0.0155  & $7.16 \times 10^{-9}$ & 4.02 & 0.0057 \\   
   0.0125 & 0.165 & $1.05 \times 10^{-8}$ & 3.85 & 0.0231 & $4.5 \times 10^{-10}$ & 4.01 & 0.0172 \\
 \hline
\end{tabular}
\caption{Convergence table showing fourth-order convergence of the new ETDRK4RDP compared with the existing ETDRK4P22 scheme for the model test-problem 3.}
\label{tab:etdrk4rdp_conv_eg3_1}
\end{table}

\begin{table}[ht]
\begin{tabular}{ p{1cm}|p{1cm}||p{2cm}|p{1cm}|p{2cm}||p{2cm}|p{1cm}|p{2cm} }
 \hline
 \multicolumn{8}{c}{ETDRK4RDP \quad \quad \quad \quad \quad \quad \quad \quad \quad ETDRK4P22} \\
 \hline
 \hline
   $k$      & $h$    & Error                 & $p$    & Cpu Time (sec) & Error           & $p$    & Cpu Time (sec) \\
   \hline
   0.10   & 0.165 & $1.5 \times 10^{-9}$ & $-$ & 0.0038   & $1.82 \times 10^{-1}$ & $-$ & 0.0015 \\       
   0.05   & 0.165 & $3.1 \times 10^{-10}$ & 2.21 & 0.0060  & $1.06 \times 10^{-2}$ & 4.09 & 0.0029 \\
   0.025  & 0.165 & $3.6 \times 10^{-11}$ & 3.11 & 0.0117  & $1.78 \times 10^{-5}$ & 9.23 & 0.0055 \\   
   0.0125 & 0.165 & $3.2 \times 10^{-12}$ & 3.50 & 0.0232 & $1.3 \times 10^{-12}$ & 23.68 & 0.0112 \\
 \hline
\end{tabular}
\caption{\label{tab etdrk4rdp_conv_eg4_1} Convergence table showing convergence of the new ETDRK4RDP compared with the existing ETDRK4P22 scheme for the model test-problem 4.}

\end{table}




\newpage
\bibliographystyle{elsarticle-num}

\begin{thebibliography}{10}

\bibitem{leveque2007finite}
R.~J. LeVeque, Finite difference methods for ordinary and partial differential equations: steady-state and time-dependent problems, SIAM, 2007.

\bibitem{tyson1999model}
R.~Tyson, S.~Lubkin, J.~D. Murray, Model and analysis of chemotactic bacterial patterns in a liquid medium, Journal of Mathematical Biology 38 (1999) 359--375.

\bibitem{kamel2007chemical}
J.~K. Kamel, Chemical vapor deposition/chemical vapor infiltration of pyrocarbon in porous carbon, University of Notre Dame, 2007.

\bibitem{zhang2022fourier}
K.~Zhang, Y.~Zuo, H.~Zhao, X.~Ma, J.~Gu, J.~Wang, Y.~Yang, C.~Yao, J.~Yao, Fourier neural operator for solving subsurface oil/water two-phase flow partial differential equation, Spe Journal 27~(03) (2022) 1815--1830.

\bibitem{anderson2000mathematical}
A.~R. Anderson, M.~A. Chaplain, E.~L. Newman, R.~J. Steele, A.~M. Thompson, Mathematical modelling of tumour invasion and metastasis, Computational and mathematical methods in medicine 2~(2) (2000) 129--154.

\bibitem{chaplain1993model}
M.~A. Chaplain, A.~M. Stuart, A model mechanism for the chemotactic response of endothelial cells to tumour angiogenesis factor, Mathematical Medicine and Biology: A Journal of the IMA 10~(3) (1993) 149--168.

\bibitem{gatenby1996reaction}
R.~A. Gatenby, E.~T. Gawlinski, A reaction-diffusion model of cancer invasion, Cancer research 56~(24) (1996) 5745--5753.

\bibitem{asante2025fourth}
E.~Asante-Asamani, A.~Kleefeld, B.~A. Wade, A fourth-order exponential time differencing scheme with dimensional splitting for non-linear reaction-diffusion systems, Journal of Computational and Applied Mathematics (2025) 116568.

\bibitem{mathews2004numerical}
J.~H. Mathews, K.~D. Fink, et~al., Numerical methods using {MATLAB}, Vol.~4, Pearson prentice hall Upper Saddle River, NJ, 2004.

\bibitem{kassam2005fourth}
A.-K. Kassam, L.~N. Trefethen, Fourth-order time-stepping for stiff {PDE}s, SIAM Journal on Scientific Computing 26~(4) (2005) 1214--1233.

\bibitem{kleefeld2012etd}
B.~Kleefeld, A.~Khaliq, B.~Wade, An {ETD} {Crank-Nicolson} method for reaction-diffusion systems, Numerical Methods for Partial Differential Equations 28~(4) (2012) 1309--1335.

\bibitem{yousuf2012numerical}
M.~Yousuf, A.~Khaliq, B.~Kleefeld, The numerical approximation of nonlinear {Black-Scholes} model for exotic path-dependent {American} options with transaction cost, International Journal of Computer Mathematics 89~(9) (2012) 1239--1254.

\bibitem{yousuf2009efficient}
M.~Yousuf, Efficient {L}-stable method for parabolic problems with application to pricing {American} options under stochastic volatility, Applied Mathematics and Computation 213~(1) (2009) 121--136.

\bibitem{khaliq2009smoothing}
A.~Khaliq, J.~Martin-Vaquero, B.~Wade, M.~Yousuf, Smoothing schemes for reaction-diffusion systems with nonsmooth data, Journal of Computational and Applied Mathematics 223~(1) (2009) 374--386.

\bibitem{asante2016real}
E.~Asante-Asamani, A.~Khaliq, B.~A. Wade, A real distinct poles exponential time differencing scheme for reaction-diffusion systems, Journal of Computational and Applied Mathematics 299 (2016) 24--34.

\bibitem{gibou2005fourth}
F.~Gibou, R.~Fedkiw, A fourth-order accurate discretization for the {Laplace} and heat equations on arbitrary domains, with applications to the {Stefan} problem, Journal of {Computational Physics} 202~(2) (2005) 577--601.

\bibitem{hundsdorfer2013numerical}
W.~Hundsdorfer, J.~G. Verwer, Numerical solution of time-dependent advection-diffusion-reaction equations, Vol.~33, Springer Science \& Business Media, 2013.

\bibitem{voss1996time}
D.~A. Voss, A.~Q.~M. Khaliq, Time-stepping algorithms for semidiscretized linear parabolic {PDEs} based on rational approximants with distinct real poles, Advances in Computational Mathematics 6 (1996) 353--363.

\bibitem{gauthier2014lectures}
P.~M. Gauthier, Lectures on several complex variables, Springer, 2014.

\bibitem{huang2021derivative}
X.~Huang, Z.~Liu, C.~Wu, Derivative and higher-order cauchy integral formula of matrix functions, Open Mathematics 19~(1) (2021) 1771--1778.

\bibitem{Davis2025suitesparse}
D.~T.~A. Davis, \href{https://github.com/DrTimothyAldenDavis/SuiteSparse/releases}{{A Suite of Sparse matrix packages}} (2025).
\newline \url{https://github.com/DrTimothyAldenDavis/SuiteSparse/releases}

\bibitem{Hanyk2014mUMFPACK}
L.~Hanyk, \href{https://geo.mff.cuni.cz/~lh/Fortran/UMFPACK/}{{UMFPACK Fortran Interface}} (2014).
\newline \url{https://geo.mff.cuni.cz/~lh/Fortran/UMFPACK/}

\bibitem{asante2020second}
E.~Asante-Asamani, A.~Kleefeld, B.~A. Wade, A second-order exponential time differencing scheme for non-linear reaction-diffusion systems with dimensional splitting, Journal of {Computational Physics} 415 (2020) 109490.

\bibitem{zegeling2004adaptive}
P.~A. Zegeling, H.~Kok, Adaptive moving mesh computations for reaction-diffusion systems, Journal of Computational and Applied Mathematics 168~(1-2) (2004) 519--528.

\bibitem{muller2015methods}
J.~M{\"u}ller, C.~Kuttler, Methods and models in mathematical biology, Lecture Notes on Mathematical Modelling in Life Sciences, Springer, Berlin (2015).

\bibitem{bhatt2015locally}
H.~P. Bhatt, A.~Q.~M. Khaliq, The locally extrapolated exponential time differencing lod scheme for multidimensional reaction--diffusion systems, Journal of Computational and Applied Mathematics 285 (2015) 256--278.

\bibitem{cherruault1990stability}
Y.~Cherruault, M.~Choubane, J.~Valleton, J.~Vincent, Stability and asymptotic behavior of a numerical solution corresponding to a diffusion-reaction equation solved by a finite difference scheme ({Crank-Nicolson}), Computers \& Mathematics with Applications 20~(11) (1990) 37--46.

\bibitem{asante2016dimensional}
E.~Asante-Asamani, B.~A. Wade, A dimensional splitting of {ETD} schemes for reaction-diffusion systems, Communications in {Computational Physics} 19~(5) (2016) 1343--1356.

\end{thebibliography}

\end{document}